\documentclass[aos,MSNbibl,nameyear,seceqn,dvips]{arximspdf}
\usepackage{dcolumn}
\usepackage{graphicx}

\doi{10.1214/13-AOS1103} 
\volume{41}
\issue{2}
\pubyear{2013}
\firstpage{897}
\lastpage{922}

\makeatletter

\newcolumntype{d}[1]{D{.}{.}{#1}}
\def\cal{\mathcal}
\newcommand{\eqref}[1]{(\ref{#1})}
\newcommand{\argmin}{\operatorname{argmin}}
\newcommand{\supp}{\operatorname{supp}}
\newcommand{\Max}{\max}
\newcommand{\Min}{\min}
\newcommand{\ds}{ }
\newtheorem{thmm}{Theorem}[section]
\newtheorem{cor}[thmm]{Corollary}
\newtheorem{lem}[thmm]{Lemma}
\newproclaim{defi}[thmm]{Definition}
\newproclaim{exam}[thmm]{Example}
\newtheorem{conjecture}[thmm]{Conjecture}

\newcommand{\er}{\mathbb{R}}
\newcommand{\al}{\alpha}
\newcommand{\eps}{\varepsilon}
\newcommand{\cA}{\mathcal{A}}
\newcommand{\cF}{\mathcal{F}}
\newcommand{\cM}{\mathcal{M}}
\newcommand{\cI}{\mathcal{I}}
\newcommand{\cP}{\mathcal{P}}
\newcommand{\cS}{\mathcal{S}}
\newcommand{\cX}{\mathcal{X}}

\newcommand{\ot}{\leftarrow}

\makeatother

\begin{document}
\begin{frontmatter}

\title{Optimal discriminating designs for several competing regression
models\thanksref{T1}}
\runtitle{Discriminating designs for several competing models}
\thankstext{T1}{Supported in part by the Collaborative
Research Center ``Statistical modeling of nonlinear dynamic processes''
(SFB 823, Teilprojekt C2) of the German Research Foundation (DFG).}

\begin{aug}
\author{\fnms{Dietrich} \snm{Braess}\ead[label=e1]{dietrich.braess@rub.de}}
\and
\author{\fnms{Holger} \snm{Dette}\corref{}\ead[label=e2]{holger.dette@rub.de}}
\runauthor{D. Braess and H. Dette}
\affiliation{Ruhr-Universit\"at Bochum}
\address{Fakult\"at f\"ur Mathematik\\
Ruhr-Universit\"at Bochum\\
44780 Bochum\\
Germany\\
\printead{e1}\\
\phantom{E-mail:\ }\printead*{e2}}
\end{aug}

\received{\smonth{5} \syear{2012}}
\revised{\smonth{12} \syear{2012}}

%
\begin{abstract}
The problem of constructing optimal discriminating designs for a class of
regression models is considered. We investigate a version of the
$T_p$-optimality criterion as introduced by Atkinson and Fedorov [\textit{Biometrika}
\textbf{62} (1975a) 289--303].
The numerical construction of optimal designs is very hard and challenging,
if the number of pairwise comparisons is larger than~2.
It is demonstrated that optimal designs with respect to this type of
criteria can be obtained by solving (nonlinear) vector-valued
approximation problems.
We use a characterization of the best approximations
to develop an efficient algorithm for the determination of the optimal
discriminating designs.
The new procedure is compared with the currently available methods in
several numerical examples, and we demonstrate that the new method can
find optimal discriminating designs in situations where the currently
available procedures fail.
\end{abstract}
%
%

\end{frontmatter}

\section{Introduction} \label{sec1}

An important problem in optimal design theory is the construction of
efficient designs for model identification
in a nonlinear relation
of the form
%
\begin{equation}
Y=\eta(x,\theta)+\varepsilon. \label{1.1}
\end{equation}
In many cases there exist several plausible models
which may be appropriate for a fit to the given data.
A typical example are dose-finding studies, where various models have
been developed for describing the dose--response relation [\citet
{pinbrebra2006}].
Some of these models, which have also been discussed by \citet
{brepinbra2005}, are listed in Table~\ref{tabex1}.
In these and similar situations the first step of the data
analysis consists of the identification of an appropriate model from a
given class of
competing regression models.

The optimal design problem for model identification
has a long history.
Early work can be found in \citet{stigler1971}, who determined 
designs for discriminating between
two nested univariate polynomials by minimizing the volume of the confidence
ellipsoid for the parameters corresponding to the extension of the
smaller model.
Several authors have worked on this approach in various other classes
of nested models
[\citet{dethal1998} or \citet{songwong1999}
among others].
%
\begin{table}
\caption{Candidate dose response models as a function of dose
$x$}\label{tabex1}
\begin{tabular*}{\textwidth}{@{\extracolsep{\fill}}lc@{}}
\hline
\multicolumn{1}{@{}l}{\textbf{Model}} & \multicolumn
{1}{c@{}}{\textbf{Full model specification}}\\
\hline
Linear & $\eta_1(x,\rho_{(1)})=60+0.56x$ \\
Quadratic & $\eta_2(x,\rho_{(2)})=60+ (7/2250) x(600-x)$ \\
E$_{\max}$ & $\eta_3(x,\rho_{(3)})=60 + 294x / (25 + x)$ \\
Logistic & $\eta_4(x,\rho_{(5)})
=49.62 + 290.51/\{1+\exp[(150 - x)/45.51]\}$ \\
\hline
\end{tabular*}
\end{table}

In a pioneering paper, \citet{atkfed1975a} proposed the $T$-optimality criterion
to construct designs for discriminating between two competing
regression models.
It provides a design such that the sum of squares
for a lack of fit test is large.
\citet{atkfed1975b} extended this approach
later 
for discriminating a selected model $\eta_1$ from
a class of other regression models,
say $ \{\eta_2, \ldots, \eta_k \}$, $k \ge2$.
This concept does not require competing nested models and has found
considerable attention
in the statistical literature; see, for example,  \citet{fedorov1981},
\citet{fedkha1986} for early and \citet{ucibog2005},
\citet{loptomtra2007}, \citeauthor{atkinson2008} (\citeyear{atkinson2008,atkinson2008b}), \citet
{Tommasi09}, \citet{wiens2009} or \citet{detmelshp2012} for
some more recent references.

In general, the problem of finding $T$-optimal designs, either
analytically or numerically, is a very hard and challenging one.
Although
\citet{atkfed1975a} indicated some arguments for the
convergence of their iterative procedure,
there is no evidence that
the convergence is sufficiently fast in cases
with more than two pairwise comparisons of regression models
such that the procedure can be used in those applications.

In the present paper we construct optimal discriminating designs
for several competing regression models where none of the models
is selected in advance to be tested against all other ones.
Let $d$ denote the number of pairwise comparison of interest.
In Section~\ref{sec2} we introduce a $T_p$-optimality criterion,
which is a weighted average of
$d$ different $T$-optimality criteria corresponding to these pairs.
It is demonstrated in Section~\ref{sec2a} that the corresponding
optimal design problems are closely
related to (nonlinear) vector-valued approximation problems.
The support points of optimal discriminating designs
are contained in the set of extreme
points of a best approximation, and
the optimal design can be determined
with the knowledge of these points.
Because we are only aware of the work of \citet{bros1968} on
vector-valued approximation,
we consider this problem in Section~\ref{sec3}.

Duality theory is then used to determine not only the points of the support,
but also the masses. The theory shows that there exist optimal designs
with a support of at most $n+1$ points,
where $n$ is the total number of
parameters in the competing regression models. We will illustrate by a
simple example that the number of support points
is usually much smaller.
It turns out that this fact occurs in particular for $d\ge3$ comparisons,
and therefore our investigations explain the difficulties in the
computation 
of $T$-optimal discriminating designs.
For this reason we find
numerical results in the literature mainly for the cases $d=1$ and $d=2$,
and advanced techniques are required
for the determination of $T_p$-optimal discriminating designs if $d \ge3$.

In Sections~\ref{sec4} and~\ref{sec5} we use the theoretical results
to develop an efficient algorithm for calculating
$T_p$-optimal discriminating designs. The main idea of the algorithm is
very simple
and essentially consists of two steps.
\begin{longlist}[(1)]
\item[(1)] The relation to the corresponding vector-valued
approximation problem is used to identify a reference set which
contains all support of the $T_p$-optimal discriminating design. This
is done by linearizing the optimization problem.
A combinatorial argument in connection with dual linear programs
determines which points are included in the support of the optimal design.
\item[(2)] A linearization of a
saddle point problem that is concealed behind the design problem
is used for a simultaneous update of all weights.
\end{longlist}
The implementation of these
two steps which are usually iterated is more complicated and described
in Section~\ref{sec5}.
Some comments regarding the convergence
and details for the main technical step of the algorithm are given in
the supplementary material [\citet{supp}] to this paper.
In Section~\ref{sec6} we provide several numerical examples and compare
our approach with the currently available methods.
In particular, we consider the problem of determining optimal
discriminating designs for the dose response models specified in Table~\ref{tabex1}.
Here the currently available procedure
fails
in the case of many pairwise comparisons, while the new method
determines a design with high efficiency
in less than $10$ iteration steps.

\section{Preliminaries} \label{sec2}

Following \citet{kiefer1974} we consider 
designs that are defined as probability measures with finite support
on a compact design space $\cX$.
If the design $\xi$ has masses $w_1,\ldots, w_\nu$
at the distinct points $x_1,\ldots, x_\nu$,
then observations are taken at these points
with the relative proportions given by the masses.
Let $\mathcal{M}= \{\eta_1,\ldots,\eta_k\}$ denote a class
of possible models for the regression function $\eta$ in \eqref{1.1},
where $\theta_{(j)}$ denotes the vector of parameters in model $\eta_j$
that varies in the set $\Theta^{(j)}$  $(j=1,\ldots,k)$.
\citet{atkfed1975b} proposed to select one model in $\mathcal{M}$, say
$\eta_1$, to fix its vector of parameters $\rho_{(1)}$ and to determine
a discriminating design by maximizing
%
\begin{equation}
\label{maximin} \min_{2\le j\le k}T_{1,j}(\xi),
\end{equation}
where
\[
T_{1,j}(\xi) :=\inf_{\theta_{(j)}\in\Theta^{(j)}} \int_\mathcal{X}
\bigl[ \eta_1 (x, \rho_{(1)})- \eta_j (x,
\theta_{(j)})\bigr]^2 \,d \xi(x) \qquad(2 \leq j \leq k).
\]
If the competing regression models $\eta_1,\ldots,\eta_k$ are not nested
(as in Table~\ref{tabex1}),
it is not clear which model should be fixed in this approach,
and it is
useful to have more ``symmetry'' in this concept.
For illustration consider
the case of two competing nonnested models, say
{$\eta_i(x,\theta_{(i)}), \eta_j(x,\theta_{(j)})$,}
and assume that the experimenter can fix a parameter for each model,
say $\rho_{(1)}$ and $\rho_{(2)}$.
In this case for a given design $\xi$ there exist two $T$-optimality
criteria, say $T_{1,2}$ and $T_{2,1}$, corresponding to the
specification of the model $\eta_1$ or $\eta_2$,
respectively, where
\begin{eqnarray}
\label{det1}
T_{i,j} (\xi) &:=& \inf_{\theta_{(j)}\in\Theta^{(j)}}
\Delta_{i,j} (\theta_{(j)},\xi)
\nonumber
\\[-8pt]
\\[-8pt]
\nonumber
&\phantom{:}= &\inf_{\theta_{(j)}\in\Theta^{(j)}}
\int_\cX\bigl[\eta_i(x,\rho_{(i)})-
\eta_j(x,\theta_{(j)})\bigr]^2 \,d \xi(x)
\end{eqnarray}
$(i\neq j)$. The first index $i$ in the term $\Delta_{i,j}$
corresponds to the fixed model $\eta_i(x,\rho_{(i)})$, while the minimum
in \eqref{det1} is taken with respect to the parameter of the model
specified by the index $j$.
The parameter associated to the minimum is denoted as
%
\begin{equation}
\label{theta-xi} \theta^*_{(i,j)} := \mathop{\argmin}_{\theta_{(j)}\in\Theta^{(j)}}
\Delta_{i,j} (\theta _{(j)},\xi),
\end{equation}
where we assume its existence and do not reflect its dependence on the
design $\xi$ and the parameter $\rho_{(i)}$
since this will always be clear from the context. {Note that we use
the notation $\theta_{(i,j)}^*$}
for the parameter corresponding to the best\vspace*{1pt} approximation of the model
$\eta_i$ (with fixed paramater $\rho_{(i)}$) by
the model~$\eta_j$.

If a discriminating design has to be constructed for $k$ competing models,
there exist $k(k-1)$ expressions of the form \eqref{det1}.
Let $p_{i,j}$ be given nonnegative weights satisfying $\sum_{i \neq j}
p_{i,j}=1$, then
a design $\xi^*$ is called \textit{$T_p$-optimal discriminating}
for the class of models $\cM= \{\eta_1,\ldots,\eta_k \}$
if it maximizes the functional
%
\begin{equation}
\label{det2} T (\xi) := \sum_{1 \leq i \neq j \leq k} p_{i,j}
T_{i,j} (\xi) = \sum_{1 \leq i \neq j \leq k} p_{i,j}
\inf_{\theta_{(j)}\in\Theta^{(j)}} \Delta_{i,j} (\theta_{(j)},\xi)
\end{equation}
[see also \citet{atkfed1975b}].
Note that the special choice $
p_{i,j} > 0$ $(j=2,\ldots,k)$, $ p_{i,j} = 0 $
$ (i=2,\ldots,k,$ $j=2,\ldots,k; i \neq j)$,\vadjust{\goodbreak}
refers to the case where one model (namely $\eta_1$) has been fixed and
is tested against all other ones.
The criterion \eqref{det2} provides a more symmetric formulation
of the general discriminating design problem. It has also been
investigated by \citet{tomlop2010} among others for $k=2$ competing
regression models.
They proposed to maximize a weighted mean of efficiencies which
is equivalent to the \mbox{criterion} (\ref{det2}) if the weights $p_{i,j}$
are chosen appropriately.

In order to deal with the general case we
denote the set of indices corresponding to the positive weights in
\eqref{det2} as
\[
\cI:=\bigl\{ (i,j) \mid p_{i,j} > 0; 1\le i \neq j \le k \bigr\}.
\]
We assume without loss of generality that the set $\mathcal{I}$ can be
decomposed in $p\leq k$ subsets
of the form $\{ (i,j) \in\cI\mid1 \le j \le k\}$ and define $\cI_i=
\{ j \in\{ 1,\ldots,k \} \mid(i,j) \in\cI\}$ as the set of indices
corresponding to those models which are used for a comparison with
model $\eta_i$.
{For each model $\eta_i$ $(i=1,2,\ldots,p)$,
a parameter, say $\rho_{(i)}$, is fixed due to prior information,}
and the model $\eta_i (x, \rho_{(i)}) $ has
to be discriminated from the other ones in the set $\cI_i$.
Define
%
\begin{equation}
\label{lambda} \lambda_i := \# \cI_i,\qquad d:=\sum
^p_{i=1} \lambda_i
\end{equation}
as the cardinality of the sets $\cI_i$ and $\cI$, respectively. Note
that $d$ denotes the total number of pairwise comparisons included in
the optimality \mbox{criterion} \eqref{det2}.
Consider the space $\cF_d= C(\mathcal{X})^d$ of continuous
vector-valued functions defined on $\cX$, and define for
a function
$ g= (g_{ij})_{(i,j)\in\cI} \in\cF_d$
a norm by
%
\begin{equation}
\label{norm} \| g \| := \sup_{x \in\cX} \bigl| g(x)\bigr |, 
\end{equation}
where $ |g(x)|^2:= \sum_{(i,j) \in\cI} p_{i,j} g^2_{ij}(x)$
denotes a \textit{weighted Euclidean norm} on~$\er^d$.
In this framework the distance between two functions $f,g \in\mathcal{F}_d$
is given by $\|f-g\|$.
{Next, given the parameters $\rho_{(1)},\ldots,\rho_{(p)}$
for the models $\eta_1,\ldots,\eta_p$, respectively, due to prior information,}
define the $d$-dimensional vector-valued function
\begin{equation}
\label{det3}
\qquad\eta(x) := \bigl(\underbrace{\eta_1(x,
\rho_{(1)}),\ldots, \eta_1(x,\rho_{(1)})}_{\lambda_1\ \mathrm{times}},
\ldots, \underbrace{\eta_p(x,\rho_{(p)}),\ldots,
\eta_p(x,\rho_{(p)})}_{\lambda_p\ \mathrm{times}} \bigr)^T,
\end{equation}
where each function $\eta_j (x,\rho_{(j)})$ appears $\lambda_j$ times
in the vector $\eta(x)$. We also
consider a vector of approximating functions
%
\begin{equation}
\label{det4} \qquad\eta(x,\theta) := \bigl( \underbrace{\bigl(\eta_j (x,
\theta_{(1,j)})\bigr) _{j \in\cI_1}}_{\in\mathbb{R}^{\lambda_1}},\ldots, \underbrace{
\bigl(\eta_j (x,\theta_{(p,j)})\bigr)_{j \in\cI_p}}_{\in\mathbb
{R}^{\lambda_p}}
\bigr)^T \in\cF_d .
\end{equation}
We emphasize again that we use the notation $\theta_{(i,j)}$
for the parameter in the model~$\eta_j$. This means that different
parameters $\theta_{(i,j)}$ and $\theta_{(k,j)}$ are used if the model
$\eta_j$ has to be discriminated from the models $\eta_i$ and $\eta_k$
($i \neq k$).
The corresponding parameters are collected in the vector
%
\begin{equation}
\label{det4a} \qquad\theta= \bigl( (\theta_{(1,j)})_{j \in\cI_1},\ldots, (
\theta_{(p,j)})_{j \in\cI_p} \bigr)^T \in\Theta =
\bigotimes^p_{i=1} \bigotimes_{j \in\cI_i} \Theta^{(j)},
\end{equation}
and we denote by
$
n := \dim\Theta= \sum^p_{i=1} \sum_{j \in\mathcal{I}_i} \dim
\Theta^{(j)}
$
the total number of all parameters involved in the $T_p$-optimal
discriminating design problem.
With this notation the optimal design problem can be rewritten as
%
\begin{equation}
\label{sattel} \max_\xi \sum_{1 \leq i \neq j \leq p}
p_{i,j} \min_{\theta_{(j)} \in\Theta^{(j)}} \Delta_{i,j} (
\theta_{(j)},\xi),
\end{equation}
and the following examples illustrate this general setting.

%
\begin{exam} \label{exam1}
Consider the case $k=3$ and
assume that all weights $p_{i,j}$ in the criterion \eqref{det2} are positive.
Here no model is preferred, and there are $6$ pairwise comparisons.
This yields $p=k=3$,
\begin{eqnarray*}
\cI&=& \bigl\{ (1,2), (1,3), (2,1), (2,3), (3,1), (3,2) \bigr\} ,
\\
\cI_1 &=& \{ 2,3 \},\qquad \cI_2 = \{ 1,3 \},\qquad
\cI_3 = \{ 1,2 \},
\end{eqnarray*}
$\lambda_1 = \lambda_2 = \lambda_3 = 2$ and $d=6$.
We obtain for the vectors in \eqref{det3} and \eqref{det4}
\begin{eqnarray*}
\eta(x) &=& \bigl( \eta_1(x,\rho_{(1)}),
\eta_1 (x, \rho_{(1)}), \eta_2 (x,
\rho_{(2)}), \eta_2 (x,\rho_{(2)}),\\
&&\hspace*{112pt}\eta_3 (x,\rho_{(3)}), \eta_3(x,
\rho_{(3)}) \bigr)^T,
\\
\eta(x,\theta) &=& \bigl( \eta_2(x,\theta_{(1,2)}),
\eta_3 (x, \theta_{(1,3)}), \eta_1 (x,
\theta_{(2,1)}), \eta_3 (x,\theta_{(2,3)}),
\\
&&\hspace*{124pt}\eta_1 (x,\theta_{(3,1)}), \eta_2(x,
\theta_{(3,2)}) \bigr)^T,
\end{eqnarray*}
with
\begin{eqnarray*}
\theta&=& (\theta_{(1,2)}, \theta_{(1,3)}, \theta_{(2,1)},
\theta_{(2,3)}, \theta_{(3,1)}, \theta_{(3,2)})^T
\\
&\in&\Theta^{(2)} \times\Theta^{(3)} \times\Theta^{(1)}
\times \Theta ^{(3)} \times\Theta^{(1)} \times
\Theta^{(2)}.
\end{eqnarray*}
\end{exam}

%
\begin{exam} \label{exam2}
Consider the problem of discriminating
between $k=3 $ nested polynomial models
$ \eta_1 (x,\theta_{(1)}) = \theta_{10} + \theta_{11}x$,
$  \eta_2 (x,\theta_{(2)}) = \theta_{20} + \theta_{21}x + \theta
_{22}x^2$ and
$\eta_3 (x,\theta_{(3)}) = \theta_{30} + \theta_{31}x + \theta_{32}x^2
+ \theta_{33}x^3$.
A common strategy to identify the degree of the polynomial is to test a
quadratic against a linear and a cubic against the quadratic model. In
this case we choose
only two positive weights $p_{2,1}$ and $p_{3,2}$ in the criterion
\eqref{det2}
which yields $
\cI= \{(2,1),(3,2) \}, $ $ \cI_1= \{ 1 \}, $ $ \cI_2 = \{ 2 \} ,$
$p=2$, $ \lambda_1=\lambda_2=1$ and $d=2$.
The functions $\eta$ and $\eta(\cdot, \theta)$ are given by $ \eta
( x)
= ( \eta_2 (x, \rho_{(2)}) , \eta_3 (x, \rho_{(3)}) )^T, $
\begin{eqnarray*}
\eta(x,\theta) &=& \bigl( \eta_1(x,\theta_{(2,1)}),
\eta_2(x,\theta _{(3,2)}) \bigr)^T = \bigl(
\theta_{10} + \theta_{11}x, \theta_{20} +
\theta_{21}x + \theta _{22}x^2
\bigr)^T,
\end{eqnarray*}
respectively, where $\theta=(\theta_{(2,1)}, \theta_{(3,2)})^T \in
\mathbb{R}^5.$
\end{exam}

\section{Characterization of optimal designs}\label{sec2a}

The $T_p$-optimality of a given design~$\xi$
can be checked by an equivalence theorem (Theorem~\ref{eqivthm})
that can be proved by the same arguments as used by \citet{atkfed1975b}.
As usual, the following properties tacitly are assumed to hold:
\begin{longlist}[(A1)]
\item[(A1)]
The regression functions $\eta_i(x,\theta_{(i)})$ are differentiable
with respect to the parameter $\theta_{(i)}$ ($i=1,\ldots, k$).
\item[(A2)] Let $\xi^*$ be a $T_p$-optimal discriminating design.
The parameter
$\theta^* = (\theta^*_{(i,j)})^T_{(i,j)\in\cI}$
defined by \eqref{theta-xi}
exists, is unique and an interior point of $\Theta$.
\end{longlist}
Both assumptions are always satisfied in linear models.
Moreover, assumption~(A1) is satisfied for many
commonly used nonlinear regression models; see \citet{sebwil1989}.
It is usually harder to check assumption (A2) because it depends on the
individual $T_p$-optimal design.
%
\begin{thmm}[(Equivalence theorem)] \label{eqivthm}
A design $\xi$ is a $T_p$-optimal discriminating design for the class
of models $\cM$
if and only if for all $x \in\cX$,
%
\begin{equation}
\label{det0} \psi(x,\xi) := \sum_{(i,j) \in\mathcal{I}}
p_{i,j} \bigl[\eta_i(x,\rho_i) -
\eta_j\bigl(x, \theta_{(i,j)}^*\bigr)\bigr]^2 \leq
T(\xi),
\end{equation}
where $\theta^*_{(i,j)}$ is defined by \eqref{theta-xi}.
Moreover, if $\xi$ is a $T_p$-optimal discriminating design,
then equality holds in \eqref{det0} for all support points of $\xi$.
\end{thmm}

The equivalence theorem asserts that there is no gap
between the solution of the max min problem \eqref{sattel}
and the corresponding min max problem.
The following result shows that the $T_p$-optimal design problem
is intimately related to a nonlinear
vector-valued approximation problem with respect to the norm \eqref{norm}.

%
\begin{thmm} \label{thm1}
Let $\eta$, $\eta(\cdot,\theta)$ and $T(\xi)$ be defined by
\eqref{det4}, \eqref{det4a} and \eqref{det2}, respectively, then
%
\begin{equation}
\label{det5} \sup_\xi T (\xi) = \inf_{\theta\in\Theta}
\bigl\| \eta- \eta(\cdot,\theta) \bigr\|^2,
\end{equation}
that is, with $\psi$ defined in the equivalence theorem
%
\begin{equation}
\label{nogap} \sup_\xi\inf_{\theta\in\Theta} \int
_{\cX} \psi(x,\xi) \,d\xi = \inf_{\theta\in\Theta} \sup
_\xi\int_{\cX} \psi(x,\xi) \,d\xi.
\end{equation}
%
If $\xi^*$ maximizes $T(\xi)$, then
the vector $\theta^* = (\theta^{*}_{(i,j)})_{(i,j) \in\cI} $
defined in \eqref{theta-xi} satisfies
%
\begin{equation}
\label{apprxprob} \bigl\| \eta(x) - \eta\bigl(x,\theta^*\bigr) \bigr\| = \inf
_{\theta\in\Theta} \bigl\| \eta(x) - \eta(x,\theta) \bigr\| = T\bigl(\xi^*\bigr).
\end{equation}
Moreover, the support of the $T_p$-optimal discriminating
design $\xi^*$ for the class $\cM$ satisfies
%
\begin{equation}
\label{det6} \supp\bigl(\xi^*\bigr) \subset\mathcal{A} := \bigl\{ x \in\cX |  \bigl|
\eta(x)-\eta\bigl(x,\theta^*\bigr) \bigr| = \bigl\| \eta- \eta\bigl(\cdot, \theta^*\bigr) \bigr\|
\bigr\} .
\end{equation}
\end{thmm}
%

\begin{pf}
We have for any design $\tilde\xi$ the relation
\[
\inf_{\theta\in\Theta} \int_{\cX} \psi(x,\xi) \,d\tilde
\xi \le\inf_{\theta\in\Theta} \sup_\xi\int
_{\cX} \psi(x,\xi) \,d\xi,
\]
and the left-hand side of \eqref{nogap} cannot be larger
than the right-hand side, that is, $T(\tilde\xi) \le\inf_{\theta
\in
\Theta} \| \eta- \eta(\cdot,\theta) \|^2$.
Since $\tilde\xi$ is an arbitrary design, the bound holds
also for $\sup_\xi T(\xi)$.
This means in terms of \eqref{det5}
$ \sup_{\xi} T(\xi) \leq\inf_{\theta\in\Theta} \|\eta-\eta
(\cdot
,\theta)\|^2$.

Now the characterization of $T_p$-optimality in
Theorem~\ref{eqivthm} and the definition of $\theta^* = (\theta
^*_{(i,j)})_{(i,j)\in\cI}$ in
Theorem~\ref{eqivthm} yield
\begin{eqnarray*}
T\bigl(\xi^*\bigr) & \leq&\inf_{\theta\in\Theta}\bigl \| \eta- \eta(\cdot,
\theta)\bigr\|^2 \leq\bigl\| \eta- \eta\bigl(\cdot,\theta^*\bigr)
\bigr\|^2
\\
&= &\sup_{x \in\cX} \sum_{(i,j) \in\cI}
p_{i,j} \bigl[ \eta_i(x,\rho_{(i)})-
\eta_j\bigl(x,\theta^*_{(i,j)}\bigr) \bigr]^2 \leq
T \bigl(\xi^*\bigr),
\end{eqnarray*}
which proves the first part of Theorem~\ref{thm1}.
The statement on the support points of $\xi^*$ follows directly from
these considerations.
\end{pf}
%

Equality \eqref{apprxprob}
means that the parameter $\theta^*$ defined in \eqref{theta-xi}
corresponds to the best approximation of the
function $\eta$ in \eqref{det3}
by functions of the form \eqref{det4}
with respect to the norm \eqref{norm}.
If this nonlinear approximation problem
has been solved,\vspace*{1pt} and the parameter
$\bar\theta= ( (\bar\theta_{(i,j)})_{j \in\cI_1}, \ldots,
(\bar\theta_{(p,j)})_{j \in\cI_p} )$ corresponds to a best
approximation, that is,
%
\begin{equation}
\label{bar-theta} \bigl\| \eta- \eta(\cdot, \bar\theta) \bigr\|^2 = \min
_{\theta\in\Theta} \bigl\| \eta- \eta(\cdot, \theta)\bigr \|^2,
\end{equation}
it follows from Theorem~\ref{thm1} that the support of the
$T_p$-optimal discriminating
design is contained in the set $\mathcal{A}$ defined in \eqref{det6}.
In linear models and in many of the commonly used nonlinear regression models
$\theta^*$ and $\bar\theta$ are uniquely determined.

%
\begin{exam} \label{exam2a}
In Example~\ref{exam1} we considered discriminating design problems
for 3 rival models $\eta_1,\eta_2,\eta_3$
and
all weights in the optimality criterion are positive.
By Theorem~\ref{thm1} the support of the $T_p$-optimal discriminating
design problem can be found by solving the nonlinear vector-valued
approximation problem
\begin{eqnarray*}
&&\inf_{\theta\in\Theta} \bigl\| \eta- \eta(\cdot,\theta)\bigr \|^2 =
\inf \biggl\{ \sup_{x \in\cX} \sum_{1 \leq i \neq j \leq3}
p_{i,j} \bigl| \eta _i(x,\rho_{(i)}) -
\eta_j(x,\theta_{(i,j)}) \bigr|^2\Big |
\\
&& \hspace*{187pt}\theta_{(i,j)} \in\Theta^{(j)};  1 \leq i \neq j \leq3
\biggr\}.
\end{eqnarray*}
\end{exam}

The following result is an approach in this framework for the
calculation of the masses of the $T_p$-optimal discriminating design.\vadjust{\goodbreak}

%
\begin{cor}
\label{weight}
Assume that a parameter $ \bar\theta$ defined in \eqref{bar-theta}
exists and is an interior point of $\Theta$, and let
$ \nabla_{\theta_{(i,j)}} $ denote the gradient of $\eta_j$ with
respect to
$\theta_{(i,j)}$.
\begin{longlist}[(a)]
\item[(a)] If a design $\xi$ is a $T_p$-optimal discriminating design for the
class $\cM$, then
%
\begin{equation}
\label{equat} \int_\mathcal{A} \bigl( \eta_i(x,
\rho_{(i)})- \eta_j(x,\bar \theta _{(i,j)}) \bigr)
\nabla_{\theta_{(i,j)}} \eta_j (x,\theta_{(i,j)})
|_{\theta_{(i,j)}= \bar\theta_{(i,j)}} \,d \xi(x)=0
\end{equation}
holds for all $(i,j) \in\cI$.

\item[(b)] Conversely, if all competing models are linear, and the design $\xi$ satisfies~\eqref{equat} such that $\supp(\xi) \subset\cA$,
then $\xi$ is a $T_p$-optimal discriminating design
for the class $\cM$.
\end{longlist}
\end{cor}

\begin{pf}
If condition \eqref{equat} is not satisfied, there is a direction
in the parameter space $\Theta$ in which the criterion decreases.
Thus \eqref{equat} is a necessary condition.
{From} Theorem~\ref{thm1} we know that the best approximation
gives rise to a $T_p$-optimal design, and it follows from
a uniqueness argument that the condition is also sufficient in this case.
\end{pf}

\section{Chebyshev approximation of $d$-variate functions} \label{sec3}
By Theorem~\ref{thm1}, a~$T_p$-optimal discriminating design
is associated to an approximation problem
in the space of
continuous $d$-variate functions on the compact design space $\cX$
where $d$ is the number of comparisons as specified by \eqref{lambda}.
This relation can be used for the computation of $T_p$-optimal designs
and
for the evaluation of the efficiency of computed designs.

In this section we will investigate these approximation
problems in more detail for the case of linear models.
We restrict the presentation to linear models because we want to emphasize
that the main difficulties already appear in linear models if $d\ge3$.
The extension
to nonlinear regression models is straightforward and will be provided
in Section~\ref{nonlin}.

The general theory here and in the previous section
provides only the information that a $T_p$-optimal discriminating
design exists
with $n+1$ or less support points where $n=\dim\Theta$.
We will demonstrate in Section~\ref{sec4.2} that the number of support points
is often much smaller than $n+1$.
This is the reason for the difficulties in the numerical construction,
even if only linear models are involved.
In contrast to other methods [see, e.g., \citet{loptomtra2007}]
the construction via the approximation
problem has the advantage that the points of the support of the
$T_p$-optimal discriminating design are directly calculated.

\subsection{Characterization of best approximations}
We will avoid double indices for vectors and vector-valued functions
throughout this section in order to avoid confusion with matrices.\vadjust{\goodbreak}
We write $\theta=( \theta_1, \theta_2, \ldots, \theta_n)^T$
instead of $(\theta_{(i,j)})^T_{(i,j)\in\cI}$ and $(f_1,\ldots,f_d)^T$
instead of the vector $\eta(x)$
defined in \eqref{det3}.
The approximation problem is considered for a given
$d$-variate function $f=(f_1,\ldots,f_d)^T \in\cF_d= C(\cX)^d$.
It is not necessary that some components of $f$ are equal
as it occurs in the function \eqref{det3}.

In the case of linear models,
equation \eqref{det4} defines an $n$-dimensional linear subspace
%
\begin{equation}
\label{vspace} V = \Biggl\{ v = \sum_{m=1}^n
\theta_m v_m \bigg| \theta=( \theta_1,
\theta_2, \ldots, \theta_n) \in\er^n \Biggr\}
\subset\cF_d ,
\end{equation}
where $v_1,v_2,\ldots,v_n \in\cF_d $ denotes a basis of $V$, and $n$ is
the dimension of the parameter space $\Theta$ in \eqref{det4a}.
Note that $f(x)$ and $v(x)$ are $d$-dimensional vectors for $x \in
\mathcal{X}$ and $v \in V$.
Theorem~\ref{thm1} relates the $T_p$-optimal discriminating design problem
to the problem of determining the best Chebyshev approximation $u^*$ of
the function $f$ by elements of the subspace $V$, that is,
\[
\bigl\|f-u^*\bigr\| = \min_{v\in V} \|f-v\|.
\]
As stated in \eqref{norm}, the norm $\|\cdot\|$ refers to the
maximum-norm on $C(\cX)^d$,
$ \|g\| := \sup_{x \in\cX} |g(x)|$,
where the weighted Euclidean norm $|\cdot|$
and the corresponding inner product in $\er^d$ are defined by
%
\begin{equation}
\label{normd} |r|^2 := \sum_{l=1}^d
p_l |r_l|^2,\qquad \langle\tilde r,r\rangle:=
\sum_{l=1}^d p_l \tilde
r_lr_l,\qquad r,\tilde r \in\er^d
\end{equation}
[here the weights $p_l$ correspond to the weights $p_{i,j}$ used in the
definition \eqref{norm}].
Because the family $V$ defined in \eqref{vspace} is a linear space,
the classical Kolmogorov criterion [see \citet{mein1967}] can be
generalized to the problem of vector-valued approximation. The result
is easily obtained from the cited literature
if products of real or complex numbers in the proof of the
classical theorem are replaced by the
Euclidean inner products of $d$-vectors.
The nonlinear character
of the procedures for determining best approximations does not matter
at this point.

%
\begin{lem}[(Kolmogorov criterion for vector-valued approximation)]\label{lem_kolm}
Let $u \in V$ and
%
\begin{equation}
\label{D} \cA:= \bigl\{ x \in\cX| \bigl|\eps(x)\bigr| =\|\eps\| \bigr\}
\end{equation}
be the set of extreme points of the error function
$
\eps:= f -u .
$
The $d$-variate function~$u$ is a best approximation to $f$ in $V$ if and only if for all $v\in V$,
%
\begin{equation}
\label{Kolmogorov} \min_{x\in\cA} \bigl\langle \eps(x), v(x)\bigr
\rangle \le0 .
\end{equation}
\end{lem}

%

Assume that $u$ is a best approximation of the function $f$.
Condition \eqref{Kolmogorov} in the Kolmogorov criterion
means that the system of inequalities
\[
\bigl\langle\eps(x), v_0(x)\bigr\rangle > 0\qquad \forall x\in\cA\vadjust{\goodbreak}
\]
is not solvable.
Let $v_1, v_2,\ldots, v_n$ be a basis of $V$.
Using the representation
%
\begin{equation}
\label{ansatz} v(x) = \sum_{m=1}^n
\al_m v_m(x)
\end{equation}
and setting $ r_m(x) := \langle\eps(x),v_m(x)\rangle$
we obtain the unsolvable system
%
\begin{equation}
\label{Rnineq} \sum_{m=1}^n
\al_m r_m(x) > 0 \qquad\forall x \in\cA
\end{equation}
for the vector $\al=(\al_1,\al_2,\ldots,\al_n)^T \in\er^n$.
The numbers $r_m(x)$ are considered as the components of a vector $r(x)$,
and by the theorem on linear inequalities [see \citet{cheney1966}, page~19] it follows
that
the system \eqref{Rnineq} is not solvable
if and only if the origin in $\er^n$ is contained in the convex hull
of the vectors
$
\{ r(x) =(r_1(x),\ldots, r_n(x))^T, x\in\cA\}.
$
By Carath\'eodory's theorem there are $\nu\le n+1$ points $x_1,\ldots
,x_\nu\in\cA$
and numbers $w_1,\ldots,w_\nu\ge0$
such that $\sum_{i=1}^\nu w_i=1$ and
%
\begin{equation}
\label{epsv} \sum^\nu_{i=1}
w_i r(x_i) = \sum_{i=1}^\nu
w_i \bigl\langle\eps (x_i),v(x_i)\bigr
\rangle = 0 \qquad\forall v \in V.
\end{equation}

%
\begin{thmm}[(Characterization theorem)]
\label{theoChar}%
Let $u \in V$ and $\cA$ be the set of extreme points
of $\eps=f-u $.
The following statements are equivalent:
\begin{longlist}[(iii)]
\item[(i)] $u $ is a best approximation to $f$ in $V$.
\item[(ii)] There exist $\nu\le n+1$ points $x_1, x_2,\ldots, x_\nu
\in\cA$ such that for all $v\in V$
%
\begin{equation}
\label{Kolmo-finite} \min_{1\le i\le\nu} \bigl\langle
\eps(x_i) , v(x_i) \bigr\rangle \le0.
\end{equation}
\item[(iii)] There exist $\nu\le n+1$ points $x_1,x_2, \ldots, x_\nu
\in\cA$
and $\nu$ weights $w_1$, $w_2, \ldots, w_\nu\ge0$,
$\sum_{i=1}^\nu w_i =1$
such that the functional
%
%
\begin{equation}
\label{functional} \ell(g) := \frac1{\|\eps\|} \sum
_{i=1}^\nu w_i \bigl\langle
\eps(x_i), g(x_i)\bigr\rangle
\end{equation}
satisfies
%
\begin{equation}
\label{propell} \ell(\eps) = \|\eps\|,\qquad \|\ell\|=1\quad \mbox{and}\quad V \subset\ker(
\ell),
\end{equation}
where $\ker(\ell) = \{ v \in V \mid\ell(v) = 0 \}$ denotes the kernel
of the linear functional~$\ell$.
\end{longlist}
%
\end{thmm}

\begin{pf}
The equivalence of (i) and (ii) follows from the Kolmogorov criterion.
To verify the equivalence with condition (iii), let $u^*$ be a best
approximation
and $\eps^*=f-u^*$.
Define the functional \eqref{functional} with the parameters
$x_i$ and $w_i$ from \eqref{epsv}.
By the Cauchy--Schwarz inequality we obtain
$\langle\eps^* (x_i),g(x_i)\rangle \le|\eps^* (x_i)| |g(x_i)|\le\|
\eps
^* \| \|g\|$
with equality if $g=\eps^* $.
Since $\sum_i w_i=1$, it follows that $\ell(g) \le\|g\|$,
again with equality if $g=\eps^* $, and the properties in \eqref{propell}
are verified.

Finally, assume that $u \in V$, and a functional with the properties
\eqref{propell} exists. We have for any $v\in V$
\[
\|f-v\| = \|\ell\| \|f-v\| \ge\ell(f-v) = \ell(f-u)+\ell(u-v) = \|f-u\| + 0,
\]
and $u $ is a best approximation.
\end{pf}

The extreme points $x_1,x_2, \ldots, x_\nu$ and the
masses $w_1,w_2, \ldots,w_\nu$ in Theorem~\ref{theoChar}
define the $T_p$-optimal discriminating design. This follows from
part (iii) of the theorem that
is closely
related to condition \eqref{equat} in Corollary~\ref{weight}.
Indeed, assume that (iii)
in the theorem is satisfied, and consider a design $\xi^*$ with weights
$w_1,w_2,\ldots, w_\nu$ at the points $x_1,x_2,\ldots, x_\nu$.
It follows for all $v \in V$ that
$ {\|\eps^* \|} \ell(v) =
\int_{\cal A} \langle f(x) -u^*(x), v(x)\rangle \,d\xi^*(x) =0 $,
and by inserting the elements $v_1,v_2,\ldots, v_n$ of the basis of $V$
we obtain precisely condition \eqref{equat}. Consequently, there exists
a $T_p$-optimal discriminating design with at most $n+1$ support points.
As we will see in Lemma~\ref{char-dual-f}, functions satisfying
only some of the properties
in Theorem~\ref{theoChar}(iii) will also play an important role.

\subsection{The number of support points---the generic case}\label{sec4.2}
By the characterization theorem there exists an optimal design
with at most $n+1$ support points.
If the number of points in the set $\cA$ equals $n+1$,
then the masses $w_1, w_2, \ldots, w_{n+1}$
of an optimal design can be calculated by the $n$ equations \eqref{equat}
together with the normalization $\int_{\cA} \,d\xi(x)=1$.
In most real-life problems, however, the number of support points
is substantially smaller than $n+1$,
and we obtain from \eqref{equat} more equations than unknown masses.
In this case the problem is ill-conditioned and the numerical
computation of the masses will be more sophisticated.
The following example illustrates the statement on the support.

%
\begin{exam} \label{examfort}
We reconsider Example~\ref{exam2}
for the polynomial regression models.
The weights
$p_{2,1}$ and $p_{3,2}$ are chosen as positive numbers.
Since all functions are polynomials, we may assume $\cX= [-1,+1]$
without loss of generality.
{A~quadratic polynomial $f_1$ is approximated by linear polynomials
in the first component,}
and a cubic polynomial $f_2$ is
approximated by quadratic polynomials in the second component.
Therefore, $V=\cP_1\times\cP_2$, where $\cP_k$ denotes
the set of polynomials of degree $\le k$.

We note that
the character of the approximation problem does not change
if we subtract a linear polynomial from $f_1$
and a quadratic polynomial from~$f_2$.
Therefore we can assume that $f(x) = (\rho_2 x^2, \rho_3 x^3)^T$.
Symmetry arguments show that the best approximating functions
will be polynomials with the same symmetry, and we obtain the
reduced approximation problem\looseness=-1
\[
\min_{\theta_1,\theta_2 \in\mathbb R} \sup_{x\in[-1,1]} \bigl(
p_{2,1} \bigl| \rho_2 x^2- \theta_1\bigr|^2
+ p_{3,2} \bigl| \rho_3 x^3 - \theta_2
x\bigr|^2 \bigr).
\]\looseness=0

We now fix the given parameters as $\rho_2=\rho_3=1$
and the weights in the $T_p$-optimality criterion as $p_{2,1}=p_{3,2} =1/2$.
The best approximation is given by
$ u^*(x) =(1/2, x)^T$,
that is, the first component is the best approximation of the univariate
function $f_1$, and the second component interpolates $f_2$
at the extreme points of $f_1-u^*_1$.
The function $\psi(x)=|f(x) -u^*(x)|^2 = (x^6-x^4+1/4)/2$
is depicted in the left part of Figure~\ref{fig1}.
The support of the $T_p$-optimal discriminating design
$\xi^*$ is a subset of the set of extreme points
$ \cA= \{-1,0,+1\}$
of the function $|f -u^*|^2$.
The linear functional $\ell$ in Theorem~\ref{theoChar}
is easily determined as
$
\ell(g) = \sqrt2 [\frac14 g(-1) - \frac12 g(0) + \frac14 g (1) ].
$
The characterization theorem, Theorem~\ref{theoChar}, yields the associated
$T_p$-optimal discriminating design
%
\begin{equation}
\label{xi-lin} \xi^*= \pmatrix{ -1 & 0 & 1 \vspace*{2pt}
\cr
\frac{1}{4} & \frac{1}{2} & \frac{1}{4}},
\end{equation}
where the first line provides the support and the second one
the associated masses.
The degeneracy is now obvious. The dimension of the set $V\subset
\mathcal{F}_2$
is $n=5$, but the solution of the corresponding approximation problem
has only $3$ extreme
points.
This degeneracy is counter intuitive.
When univariate functions are approximated by polynomials in $\cP_2$, then
by Chebyshev's theorem there are at least 4 extreme points.
Although our approximation problem with 2-variate functions
contains more functions and more parameters,
the number of extreme points is smaller.

%
\begin{figure}

\includegraphics{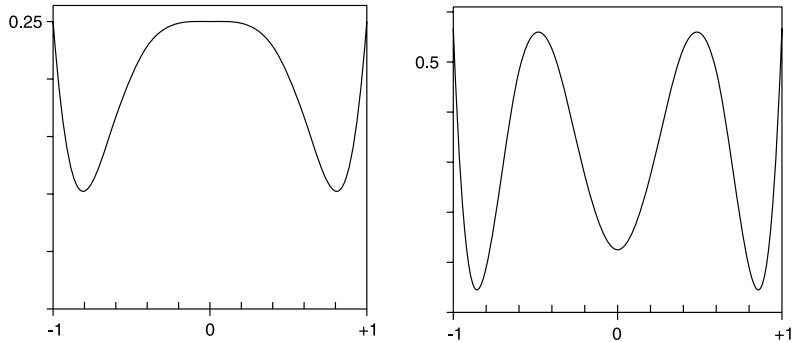}

\caption{Error functions $\psi(x)=|f(x) - u^*(x)|^2$ in the
equivalence theorem
for Example~\protect\ref{examfort}. Left panel: $\rho_2=\rho_3=1$;
right panel: $\rho_2=1, \rho_3=4$.}\label{fig1}
\end{figure}
%

Note also that
the second component is determined by interpolation
and not by a direct optimization.
The same designs are obtained whenever
$
p_{2,1} \rho_2^2 \ge p_{3,2} \rho_3^2$.
If this condition does not hold, we may have $4$ extreme points,
as shown in the right
part of Figure~\ref{fig1} for the choice $\rho_2= 1; $ $\rho_3=4$.
The solution is also degenerate. Here,
the location of the support points depends on the value of $\rho_3$.
In the mentioned case we obtain (subject to rounding)
the $T_p$-optimal discriminating design with masses
$ 0.18$, $ 0.32$, $ 0.32$, $ 0.18$
at the points $-1 $, $ -0.48$, $0.48$ and $ 1 $.
\end{exam}

The previous example shows that the cardinality of the support
depends on the given parameters
$\rho_{(1)}, \rho_{(2)}, \ldots, \rho_{(p)} $.
The following definition helps one to understand which cardinality
is found in most cases.

\begin{defi}
Let $\xi^*$ be a $T_p$-optimal discriminating design for the given data
$\rho_{(1)}, \rho_{(2)}, \ldots, \rho_{(p)} $ with
$\nu\le n+1$ support points.
The design $\xi^*$ is called a \textit{generic point}
if for all parameters in some neighborhood of
$\rho_{(1)}, \rho_{(2)}, \ldots, \rho_{(p)} $
the corresponding $T_p$-optimal discriminating designs have the same
number $\nu$
of support points.
\end{defi}
%
Our numerical experience leads to the following:

\begin{conjecture}
\label{con45}
If a $T_p$-optimal discriminating design is a generic point, then
its support consists of
\[
\max_{(i,j) \in\cI} \dim\Theta^{(j)} + 1
\]
points.
\end{conjecture}
It has been observed in the literature that the number of points in the support
can be smaller than $n+1$ [see, e.g., \citet{dettit2009}], but
computations for $d \le2$ do not
give the correct impression how large the reduction can be.

\section{Linearization and duality}
\label{sec4}
The equivalence theorem (Theorem~\ref{eqivthm}) and Theorem~\ref{thm1} show
that the maximization
of $T(\xi)$ is related to a minimization problem.
This duality is also reflected in the characterization theorem (Theorem~\ref
{theoChar}).
We will now
consider Newton's iteration for the computation
of best approximations.

In each step of the iteration an approximating function $u$
in the family $V$ is improved
simultaneously with a \textit{reference set} $\cS=\{x_1,x_2,\ldots
,x_\nu\}$
that is considered as an approximation of the set $\cA$ of extreme
points which contains the support of $T_p$-optimal discriminating designs.
Thus we focus on the minimization problem, but we will obtain the
associated weights $\{w_1,w_2,\ldots,w_\nu\}$ by duality considerations.
Note that in this section we regard duality in connection
with the linearized problems and the involved linear programs.

Given a guess $u$ for the approximating function
and a finite reference set $\cS$, the quadratic term of the correction $v$
in the binomial formula is temporarily ignored.
As usual, let $\eps:=f-u$.
We replace the optimization problem
%
\begin{eqnarray}
\label{full} &&\Max_{x_i \in\cS} \bigl|f(x_i)-u (x_i)-v(x_i)
\bigr|^2
\nonumber
\\[-8pt]
\\[-8pt]
\nonumber
&&\qquad=\Max_{x_i \in\cS} \bigl\{\bigl|\eps(x_i)\bigr|^2 - 2 \bigl
\langle\eps(x_i), v(x_i)\bigr\rangle
+\bigl|v(x_i)\bigr|^2 \bigr\} \to\Min_{v\in V}!
\end{eqnarray}
by the linear program
%
\begin{equation}
\label{LP} \max_{x_i \in\cS} \bigl\{\bigl|\eps(x_i)\bigr|^2
- 2 \bigl\langle\eps(x_i), v(x_i)\bigr\rangle \bigr\}
\to\min_{v\in V}!.
\end{equation}
While the left-hand side of \eqref{full} is obviously
bounded from below, this is not always true for the optimization
problem \eqref{LP}.
The boundedness, however, is essential for the algorithm.

%
\begin{defi}
A function $u \in V$ is called \textit{dual feasible} for the reference
set $\mathcal{S}$,
if the left-hand side of \eqref{LP} is bounded from below.
\end{defi}

The notation of dual feasibility will be clear from the dual linear
program \eqref{adjointf}
and Lemma~\ref{prim-dual} below.
We will also see in Lemma~\ref{char-dual-f} that only the dual feasible
functions
are associated to a design $\xi$ in the sense of \eqref{det1}.

The minimization of a linearized functional on a
finite set $\cS=\{x_i\}_{i=1}^\nu$ with $\nu\ge n+1$
as in \eqref{LP} will be the basis of our algorithm.
%
For a given error function $\eps=f-u$
and a reference set with $\nu$ points $x_1,x_2,\ldots,x_\nu$ we may
use representation \eqref{ansatz} and
rewrite the primal problem \eqref{LP} as
a linear program for the
$n+1$ variables $E, \al_1, \al_2,\ldots, \al_n$:
%
\begin{eqnarray}
\label{primal} %
E &\to& \min!,
\nonumber
\\[-9pt]
\\[-9pt]
\nonumber
\ds2 \sum_{m=1}^n \al_m
\bigl\langle\eps(x_i), v_m(x_i)\bigr\rangle
+ E &\ge& \bigl|\eps(x_i)\bigr|^2 , \qquad i=1,2,\ldots,\nu.
\end{eqnarray}
Obviously, there exists a feasible point for this linear program,
since the inequalities are satisfied
by $\al_1=\al_2=\cdots=\al_n=0$ and $E=\|\eps\|^2$.

The dual program to \eqref{primal} contains
the equations for the $\nu$ weights $w_1$, $w_2,\ldots, w_\nu$
with the adjoint matrix, where we can drop the factor $2$ for the sake
of simplicity,
%
\begin{eqnarray}
\label{adjointf} %
 \ds\sum
_{i=0}^\nu w_i \bigl|\eps(x_i)\bigr|^2
&\to& \max!,
\nonumber\\[-2pt]
\ds\sum_{i=1}^\nu w_i \bigl
\langle\eps(x_i), v(x_i)\bigr\rangle &=& 0\qquad \forall v
\in V,
\\[-2pt]
\ds\sum_{i=1}^\nu w_i &=& 1,\qquad
w_i \ge0, i=1,2,\ldots,\nu. \nonumber
\end{eqnarray}
The following result of duality theory will play
an important role [for a proof see \citet{pastei1998}].

%
\begin{lem}
\label{prim-dual}
The linear program \eqref{adjointf} has a feasible point
and a solution
if and only if the objective function in the linear program \eqref{primal}
is bounded from below, that is,
\[
\min_{v\in V} \max_{0\le i\le n} \bigl\langle
\eps(x_i), v(x_i)\bigr\rangle \ge0.\vadjust{\goodbreak}
\]
%
\end{lem}
%

If the linear program \eqref{adjointf} has a feasible point,
there is a solution with at most $n+1$ positive weights.
We obtain a linear functional $\ell$ of the form \eqref{functional}
with these parameters where
$ \|\ell\|=1$ and $V \subset\ker(\ell)$.
We have $\ell(\eps) < \|\eps\|$, whenever $u$ is not a best approximation.
Since the values of the primal program \eqref{primal}
and the dual program \eqref{adjointf}
coincide, we also have
\[
E = \sum_{i=1}^\nu w_i \bigl|
\eps(x_i)\bigr|^2.
\]

The final result of this section shows that the evaluation of the
functional $T$ defined in \eqref{det2} for a given
design $\xi$ is strongly related to dual feasibility.

%
\begin{lem}
\label{char-dual-f}
Let $u\in V$ and $\mathcal{S}=\{x_i\}_{i=1}^\nu$.
The following statements are equivalent:
\begin{longlist}[(iii)]
\item[(i)] The function $u$ is dual feasible for the reference set
$\mathcal{S}$.
\item[(ii)] There exist nonnegative weights $w_i, i=1,2,\ldots,\nu$,
such that
\[
\sum_{i=1}^\nu w_i \bigl
\langle\eps(x_i), v(x_i)\bigr\rangle = 0
\]
{holds} for all $ v \in V$.
\item[(iii)] There exists a design $\xi$ supported on $\mathcal{S}$
such that
\[
u = \mathop{\argmin}_{v \in V} \int_{\cX} \bigl|f (x) -v(x)
\bigr|^2 \,d\xi(x).
\]
\end{longlist}
\end{lem}


\begin{pf}
The equivalence of (i) and (ii) is a direct consequence of Lem\-ma~5.2.
Note that for $t \in\er$ and $v \in V$,
%
\begin{eqnarray}
\label{lower} &&\sum^\nu_{i=1}
w_i \bigl|(f - u - tv) (x_i)\bigr|^2
\nonumber
\\[-8pt]
\\[-8pt]
\nonumber
&&\qquad = \sum
^\nu_{i=1} w_i \bigl( \bigl|
\eps(x_i)\bigr|^2 - 2t\bigl\langle\eps(x_i),
v(x_i)\bigr\rangle + t^2 \bigl|v(x_i)\bigr|^2
\bigr).
\end{eqnarray}
If (ii) holds with the weights $w_i$, then expression
\eqref{lower}
attains its minimum at $v=0$. Hence, $u$ is the solution of the
minimization of $ \int_{\cX} |f (x) -u(x) |^2 \,d\xi(x) $ for the
design $\xi$
with the support $\mathcal{S}$ and the masses $w_1 , \ldots, w_\nu$
from condition~(ii).
If (ii) does not hold, then the minimum of \eqref{lower}
is not obtained at $t=0$ for one $v\in V$.
Therefore the minimum is not attained at $u$.
\end{pf}

\section{The algorithm}
\label{sec5}
Each step of our iterative procedure consists of two parts.
The first part deals with the improvement of the approximating
function and the reference set. It focuses on the approximation problem.
The second part is concerned with the computation of the
associated masses.
The dual linear program is embedded in\vadjust{\goodbreak} a saddle point problem.
Thus computations for the primal problem
and the dual problem may alternate during the iteration.
The small number of support points of $T_p$-optimal discriminating
designs (as described
in Conjecture~\ref{con45}) has impact on both parts.

The iteration starts with a set of parameters
$\theta_{(i,j)} $ and a reference set of about $n+1$ points
which divide the interval $\cX$ into subdomains of equal size.
Of course, any prior information may be used
for getting a better initial guess.

\subsection{Newton's method and its adaptation} \label{sec61}
The improvement of the approximation on a given reference set
will be done iteratively by Newton's method.
In order to avoid the introduction of an additional symbol,
we focus on one step of the iteration for the given
input $u_0$, the corresponding error function $\eps_0=f-u_0$,
and the reference set $\cS_0$.
The simplest Newton step,
\begin{itemize}
\item[]
\textit{Given $u_0$ and $\cS_0$, find a solution of the linear program
\textup{\eqref{primal}} for $u=u_0$, set $v = \sum_m \al_m v_m$.}

\item[]\textit{Take $u_1 = u_0+v$ as the result of the Newton step,}
\end{itemize}
looks natural; however, it can be only the basis of our algorithm.
We take three actions.
For convenience, we use the notation
$\|g\|_{\cS} := \sup_{x\in\cS} |g(x)|$.
%
%
\begin{longlist}[(1)]
\item[(1)] \textit{Newton steps on subspaces.}
Referring to the notation in Section~\ref{sec2}
we write the space of approximating functions as a sum
of $d$ subspaces
%
\begin{equation}
\label{Vij} V = \bigoplus_{(i,j)\in\cI} V_{(i,j)} ,
\end{equation}
where $V_{(i,j)}$ contains those functions in $V$ that correspond to
$\{\eta_j(\cdot,\theta_{(i,j)} ) \mid\theta_{(i,j)} \in\Theta
^{(j)}\}$.
The linear program that is obtained from \eqref{primal}
by the restriction of the functions $\sum_m \al_m v_m$
to the subspace $V_{(i,j)}$  will be denoted as \eqref{primal}$_{(i,j)}$.

The improvement of the approximation on the reference set will be
done iteratively by Newton's method. The linearization \eqref{LP},
however, will be considered for the subspaces $ V_{(i,j)}$
and not for $V$. In other words, the $d$ linear programs \eqref
{primal}$_{(i,j)}$
are performed separately.
It follows from Conjecture~\ref{con45}
that we have dual feasibility only on lower dimensional spaces.
Indeed, the splitting \eqref{Vij} creates dual feasible problems,
or the defect is one-dimensional, and the regularization
described in item (3) below is the correct remedy.
Moreover, another improvement without the splitting
will be provided in combination with the evaluation
of the masses in part 2 of the iteration step.
[Note that we have the same splitting in the evaluation
of $\theta^*$ according to \eqref{theta-xi}.]

\item[(2)] \textit{The damped Newton method.}
The Newton correction $v$ will be multiplied by a damping factor $t$.
By definition of the Newton method we have
$\max_{x_i \in\cS_0} \{ |\eps_0(x_i)|^2 - 2\langle\eps_0(x_i),
v(x_i)\rangle\}
< \|\eps_0\|_{\cS_0}^2
$
if\vadjust{\goodbreak} we have not yet obtained the solution of the actual minimum problem.
Since
\[
\bigl|(f-u_0-tv) (x_i)\bigr|^2 = \bigl|
\eps_0(x_i)\bigr|^2 - 2t \bigl\langle
\eps_0(x_i), v(x_i)\bigr\rangle + O
\bigl(t^2\bigr),
\]
it follows that $\|f-u_0-tv)\|_{\cS_0}^2 < \|\eps_0\|_{\cS_0}^2$
for sufficiently small positive factors~$t$; and thus an improvement is
generated.
Let
\[
T:=\bigl\{ 1, 2^{-1}, 2^{-2}, 2^{-3},
2^{-4}, \ldots, 2^{-7},0\bigr\},
\]
and determine
%
\begin{equation}
\label{min-t} t = \mathop{\argmin}_{t\in T} \|f-u_0-tv
\|_{\cS_0} .
\end{equation}
The standard set of damping factors $1, 2^{-1}, 2^{-2},\ldots$
has been augmented by the element $0$, and therefore
the new approximation is at least as good as the old one.

\item[(3)] \textit{Regularization by adding a bound.}
By definition the objective function $E$ is not bounded from below
in the linear program \eqref{primal}$_{(i,j)}$ if $u_0$
is not feasible with respect to $ V_{(i,j)}$.
Therefore, we add the restriction
$
E \ge0
$
to the linear programs.
\end{longlist}
%
%
At the end of this part of the iteration step
we have an improved approximation~$u_1$.
Extreme points of $f-u_1$ that are not yet obtained in $\cS_0$
are added to this set.
A decision on the augmentation of the reference set is easy
when the error curve is shown on the monitory of the computer.
Furthermore, we mark the points in $\cS_0$
to which a positive mass was given by the dual linear
program associated to \eqref{primal}$_{(i,j)}$ for one pair $(i,j) \in
\cI$.
The points in the reference set are relabel
such that $x_1,x_2,\ldots,x_\mu$ are the marked ones.

\subsection{Computation of best designs}
The adapted Newton step in the first part of the iteration step
has provided an improved error curve $\eps_1=f-u_1$
and simultaneously a set of marked points, say $\{x_1,x_2,\ldots,x_\mu
\}$.
Let $\xi$ be a design with this support
and masses $\{w_1,w_2,\ldots,w_\mu\}$ that are not yet known.
We look for a correction $v$ with the representation
\eqref{ansatz} such that $u_1+v$ is associated to $T(\xi)$
in the spirit of \eqref{det1}, that is,
we have to minimize
%
\begin{equation}
\label{sat-1} \ds\sum^\mu_{i=1}
w_i \Biggl| \eps_1(x_i) - \sum
^n_{k=1} \al_k v_k(x_i)
\Biggr|^2  = \al^T A \al-2 w^T R\al+
b^T w, 
\end{equation}
where the elements of the matrices $A=(A_{jk})_{j,k=1,\ldots,n}$,
$R=(R_{ik})^{j=1,\ldots,n}_{i=1,\ldots,\mu}$
and the vector $b=(b_i,\ldots,b_n)^T$ are defined by
%
\begin{eqnarray}\label{sat-2}
A_{jk} &:=& \sum_i w_i
\bigl\langle v_j(x_i), v_k(x_i)
\bigr\rangle,\nonumber
\\
 R_{ik} &:= &\bigl\langle \eps_1(x_i),
v_k(x_i)\bigr\rangle,
\\
b_i &:=& \bigl|\eps(x_i)\bigr|^2.
\nonumber
\end{eqnarray}
The optimal design among all designs supported at $\{ x_1,\ldots,x_\mu
\}
$ is determined
by the solution of the saddle point problem
%
\begin{equation}
\label{sat-3} \max_{w} \min_\al\bigl\{
\al' A \al-2 w'R\al+ b'w \bigr\},
\end{equation}
where we will ignore the dependence
of the matrix $A$ on $w$ for a moment.
Reasonable weights in~\eqref{sat-2} will be specified below.
The inner optimization problem in~\eqref{sat-3} is solved by
%
\begin{equation}
\label{sat-7} A\al= R^T w,
\end{equation}
and we arrive at the quadratic program
%
\begin{equation}
\label{quadprog} \max_w \bigl\{ - w^T
RA^{-1}R^T w+ b^T w | e^T w=1,
w_i \ge0\bigr\},
\end{equation}
where $e:=(1,1,\ldots,1)^T$ is a $\mu$-vector.
In order to check whether all masses are positive,
we compute an approximate solution $\tilde w$
by solving the linear program
%
\begin{eqnarray}
\label{sat-4} %
\sum
_{i=1}^\nu\bigl| \bigl(RA^{-1}R^T
\tilde w\bigr)_i \bigr| &\to& \min!,
\nonumber
\\[-8pt]
\\[-8pt]
\nonumber
\ds\sum_{i=1}^\nu\tilde w_i
&=& 1,\qquad \tilde w_i \ge0,  i=1,2,\ldots,\mu.
\end{eqnarray}

We observed in our numerical calculations that all masses are positive,
whenever at least 2 points have been marked in part 1 of the procedure.
After removing points with zero mass $\tilde w_i$, if necessary,
we can ignore the restrictions $w_i \ge0$,
and problem \eqref{quadprog}
is solved by the linear saddle point equation
%
\begin{equation}
\label{sat-5} \pmatrix{ A & -R^T & \vspace*{2pt}
\cr
-R & & e
\vspace*{2pt}
\cr
& e^T &} \pmatrix{ \al\vspace*{2pt}
\cr
w
\vspace*{2pt}
\cr
\lambda } = \pmatrix{ 0 \vspace*{2pt}
\cr
-\frac12 b \vspace*{2pt}
\cr
1}.
\end{equation}
Now we are in a position to specify which masses
are inserted in \eqref{sat-2} when the matrix $A$ is calculated.
We start with equal masses $w_i =1/\mu$ for $i=1,2,\ldots, \mu$
when we build the matrix for the linear program \eqref{sat-4}.
The masses $\tilde w_i$ from the linear program are then used
in the definition of the matrix $A$ for the saddle point equation
\eqref{sat-5}. The solution of \eqref{sat-5} yields the masses
for the improved design $\xi$.
By definition, these masses are used when the criterion $T(\xi)$ is evaluated.

The evaluation of $T(\xi)$ according to \eqref{det2}
provides also corrections of the parameters.
Let $u_2 = u_1+v$ be the associated function in $V$.
By definition the sum of weighted squares $\sum_i |\eps(x_i)|^2$
is smaller for $u_2$ than for $u_1$.
If the errors are nearly equilibrated, it follows
that $\max_i |\eps(x_i)|$ will also be smaller for $u_2$ than for $u_1$.
Therefore, we look for a damping factor $t$ such that
the norm of the error\vadjust{\goodbreak} $\|f-(u_1+tv)\|$ is as small as possible.
The details of the damping procedure are the same
as in the damped Newton method described in Section~\ref{sec61}.

The value of $T(\xi)$
is a lower bound for the degree of approximation
and
provides a lower bound of the $T_p$-efficiency
%
\begin{equation}
\label{eff1} \operatorname{Eff}_{T_p} (\xi) := \frac{T (\xi)}{\sup_\eta T (\eta)} \ge
\frac{T (\xi)} {\| \eps\|^2} .
\end{equation}
In particular, we have a stopping criterion for the algorithm. The
iteration will
be stopped if the guaranteed $T_p$-efficiency
is sufficiently close to 1.

\subsection{Adaptation to nonlinear models} \label{nonlin}

When the models $\eta_1, \eta_2,\ldots, \eta_k$ depend
nonlinearly on the parameters,
the approximating function $u(x,\theta)$ depends in a (possibly)
nonlinear way on the parameter $\theta$.
The gradient space defined by
%
\begin{equation}
\label{crit} \Biggl\{ \eta(\cdot,\tilde\theta) + \sum
_{k=1}^n \al_k \frac{\partial}{\partial\theta_k} \eta(
\cdot,\theta)\bigg|_{\theta=\tilde\theta}, \al\in\er^n \Biggr\}
\end{equation}
is a linear subspace and all the procedures described
for linear spaces can be applied to this
gradient space.
Only the computation of $T(\xi)$ for given $\xi$ requires more
effort. The minimization in its definition of $T(\xi)$ can be
done by Newton's method. The linearization uses those
formulas that are related to the minimization in the
gradient space.
Thus the algorithm can also deal with nonlinear models.

\section{Numerical results} \label{sec6}

We confirm the efficiency of the new algorithm by numerical results
for three examples with linear and nonlinear regression functions. A~fourth
example can be found in Appendix~C of the supplementary material [\citet{supp}].
We also provide a comparison with the algorithm proposed by \citet{atkfed1975a}.
Each iteration step is performed in the examples
in less than 1 or 2 seconds on a five years old personal computer.
The quotient $T(\xi_j)/\|\eps_j\|^2$ in the tables shows the lower
bound for
the efficiency defined in~\eqref{eff1}.
When we distinguish between part 1 and part~2 of the iteration step,
an index is added to the iteration count. In particular, we distiguish
the error functions $\eps_{j,1}$
and $\eps_{j,2}$ obtained in part 1 and part~2 of the iteration.
The ratio $T(\xi_j)/\|\eps_{j,2}\|^2$ in the tables shows the lower
bound for the efficiency defined in~\eqref{eff1}.

%
\begin{table}
\tabcolsep=0pt
\caption{The results of the new algorithm for Example \protect\ref
{exam6.1}}\label{table61}
\begin{tabular*}{\textwidth}{@{\extracolsep{4in minus 4in}}ld{2.4}ccccc@{}}
\hline
& \multicolumn{1}{c}{\textbf{Part 1}} & \multicolumn{3}{c}{\textbf
{Part 2}} & & \\[-4pt]
& \multicolumn{1}{c}{\hrulefill} & \multicolumn{3}{c}{\hrulefill} &
& \\
\multicolumn{1}{@{}l}{$\bolds{j}$} & \multicolumn{1}{c}{$\bolds{\|
\eps_{j,1}\|^2}$} &
\multicolumn{1}{c}{$\bolds{\|\eps_{j,2}\|^2}$} & \multicolumn
{1}{c}{$\bolds{T(\xi_j)}$}
& \multicolumn{1}{c}{$\bolds{\ds\frac{T (\xi_j)}{\|\eps_{j,2} \|^2}}$}
& \multicolumn{1}{c}{\textbf{Support}} & \multicolumn
{1}{c@{}}{\hspace*{28pt}\textbf{Reference set}}\\
\hline
$0$ & 12.5 & & & & &$\cS= \{-1,-0.5,$
\\
& & & & & &\hspace*{56pt}$ -0.1,0,0.1,0.5,1\}$
\\
$1$ & 2.3513 & & & & & \\
$2$ & 0.6092 & & & & & \\
$3$ & 0.3391 & 0.2434 & 0.0146 & 0.0600 & $\{-1,0.3,1\}$
& $\cS\ot\cS\cup\{0.3\}$\\
$4$ & 0.2012 & 0.1556 & 0.1144 & 0.7350 & $\{-1,0.2,1\}$
& $\cS\ot\cS\cup\{0.2\}$ \\
$5$ & 0.1401 & 0.1287 & 0.1029 & 0.8002 & $\{-1,-0.3,1\}$
& $\cS\ot\cS\cup\{-0.3\}$ \\
$6$ & 0.1268 & 0.1265 & 0.1225 & 0.9685 & $\{-1,0.1,1\}$
& \\
$7$ & 0.1261 & 0.1260 & 0.1244 & 0.9872 & $\{-1,-0.05,1\}$
& $\cS\ot\cS\cup\{-0.05\}$\\
$8$ & 0.1259 & 0.1258 & 0.1246 & 0.9906 & $\{-1,0,04\}$
& $\cS\ot\cS\cup\{0.04\}$ \\
\hline
\end{tabular*}
\end{table}

%

%
\begin{exam}
\label{exam6.1}
We consider once more Example~\ref{exam2},
fix $p_{2,1}=p_{3,2}= \frac12$, set
\[
f(x) = \bigl(\eta_2(x,\rho_{(2)}), \eta_3(x,
\rho_{(3)})\bigr)^T = \bigl(1+x+x^2,
1+x+x^2+x^3\bigr)^T
\]
and start the algorithm with $u_0 = (0, 0)^T$,
that is, $\theta_{(2,1)} = (0, 0)$, $\theta_{(3,2)} = (0, 0, 0)$.
The initial guess $u_0$ implies that the functions obtained during the iteration
do not have the symmetry properties discussed in Example~\ref{examfort}.\vadjust{\goodbreak}

The results of the new algorithm are displayed in Table~\ref{table61}.
After $8$ iteration steps we obtain a
discriminating design with at least 99\% efficiency. In the first part
of the iteration
the lower bound is very small and of no use, but it is increasing
rapidly during the iteration.
In Figure~\ref{fig2} we display the shape of the error function in the
first $3$ iterations. We observe that
the location of the extreme points changes substantially
in the first steps of the algorithm. A comparison with Figure \ref
{fig1} shows
that afterwards there are no substantial changes of the shape.
The resulting discriminating design
puts the masses $0.241$, $0.501$ and $0.258$
at the points $-1$, $0.04$ and~$1$, respectively.
The parameters may be compared with the exact optimal ones $ \frac
{1}{4} $, $\frac{1}{2} $ and $ \frac{1}{4}$ at the points given in
\eqref{xi-lin}.
The parameters corresponding to the solution of the nonlinear
approximation problem defined by the right-hand side of \eqref{det5}
are given by $\overline\theta_{(2,1)} = (1.501, 1.002)$,
$\overline\theta_{(3,2)} = (0.996, 1.976, 0.958)$.

\begin{figure}[b]

\includegraphics{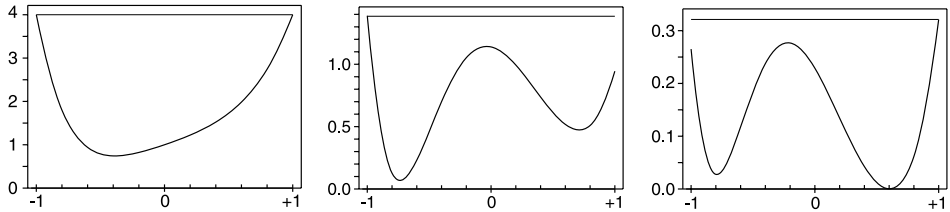}

\caption{Error curve $|f - u|^2$ in the first three iteration steps for
Example \protect\ref{exam6.1}.}\label{fig2}
\end{figure}

%
\begin{table}
\caption{The results of the algorithm proposed by Atkinson and Fedorov
(\citeyear{atkfed1975a})
in Example \protect\ref{exam6.1} (left part) and Example \protect
\ref
{exam6.2} (right part)}\label{table61alt}
\begin{tabular*}{\textwidth}{@{\extracolsep{\fill}}l@{\hspace*{40pt}}ccc@{\hspace*{40pt}}ccc@{}}
\hline
$\bolds{j} $ & $\bolds{\Vert \psi\Vert }$ & $\bolds{T(\xi_j)}$ & $\bolds
{\ds\frac{T(\xi_j)}{\Vert \psi\Vert }}$&
$\bolds{\Vert \psi\Vert }$ & $\bolds{\ds T(\xi_j)}$ & $\bolds{\ds\frac
{T(\xi_j)}{\Vert \psi\Vert }}$\\
\hline
\phantom{0}1 & 0.2172 & 0.1041 & 0.4791 &0.0104 &0.0033 &0.3150 \\
\phantom{0}2 & 0.3995 & 0.0743 & 0.1860 &0.0133 &0.0034 &0.2560\\
\phantom{0}3 & 0.3189 & 0.0778 & 0.2440 &0.0241 &0.0045 &0.1880\\
\phantom{0}4 & 0.1539 & 0.1216 & 0.7903 &0.0099 &0.0055 &0.5583 \\
\phantom{0}5 & 0.1974 & 0.1195 & 0.6055 &0.0131 &0.0055 &0.4199\\
\phantom{0}6 & 0.2337 & 0.1137 & 0.4868 &0.0094 &0.0063 &0.6682\\
\phantom{0}7 & 0.2045 & 0.1143 & 0.5592 &0.0093 &0.0060 &0.6471\\
\phantom{0}8 & 0.1347 & 0.1240 & 0.9206 &0.0121 &0.0062 &0.5161\\
\phantom{0}9 & 0.1732 & 0.1217 & 0.7029 &0.0079 &0.0065 &0.8228\\
10 & 0.2055 & 0.1186 & 0.5773 & 0.0104 &0.0064 &0.6153\\
11 & 0.1791 & 0.1199 & 0.6694 &0.0099 &0.0064 &0.6502\\
12 & 0.1356 & 0.1243 & 0.9165 &0.0091 &0.0065 &0.7166\\
13 & 0.1651 & 0.1200 & 0.7267 &0.0081 &0.0066 &0.8223\\
14 & 0.1640 & 0.1229 & 0.7493 &0.0088 &0.0065 &0.7371\\
15 & 0.1714 & 0.1213 & 0.7078 &0.0097 &0.0065 &0.6686\\
16 & 0.1362 & 0.1243 & 0.9130 &0.0070 &0.0067 &0.9550\\
\hline
\end{tabular*}
\end{table}

%

For the sake of comparison we also present in the left part of Table~\ref{table61alt}
the corresponding results for the first $16$ iterations of the
algorithm proposed
by \citet{atkfed1975a}.
This method starts with an initial guess, say $\xi_0$, and computes
successively new designs $\xi_1, \xi_2, \ldots$ as follows:

%
%
\begin{longlist}[(1)]
\item[(1)] At stage $s$ a point $x_{s+1} \in\mathcal{X}$ is determined
such that
$
\psi(x_{s+1},\xi_s) = \sup_{x \in\mathcal{X}} \psi(x, \xi_s),
$
where the function $\psi$ is defined in \eqref{det0}.
\item[(2)] The updated design $\xi_{s+1}$ is defined by
$
\xi_{s+1} = (1- \alpha_s)\xi_s + \alpha_s \delta_{x_{s+1}},
$
where $\delta_x$ is the Dirac measure at point $x$, and $(\alpha_s)_{s
\in\mathbb{N}_0}$ is any sequence of positive numbers satisfying
$
\alpha_s \rightarrow0; $
$ \sum^\infty_{s=0} \alpha_s = \infty; $ $ \sum^\infty_{s=0}
\alpha
^2_s < \infty.$
\end{longlist}
%
This procedure provides the design
\[
\xi= \pmatrix{ -1 & -0.2 & -0.1 & 0 & 0.1 & 1\vspace*{2pt}
\cr
0.23 & 0.18 & 0.12
& 0.1 & 0.17 & 0.20 }
\]
in $12$ iteration steps, and its efficiency is at least $92\%$.
The final design contains an unnecessarily large support,
although several design points with low weight have been removed during
the computations.
Note that neither the sup-norm of the function
$\psi$ is decreasing, nor the lower bound $T(\xi)/\| \psi\|$ is
increasing. In particular if the iteration is continued,
the lower bound for the efficiency of the calculated design is
decreasing again. This effect is at first compensated
after the~$16$th iteration, where the bound for the efficiency is $91\%
$ (but not $92\%$ as
after the $12$th iteration). This ``oscillating behavior'' was also
observed in other examples and seems to be typical for the
{frequently used} algorithm proposed by \mbox{\citet{atkfed1975a}}.
\end{exam}

%
\begin{table}
\tabcolsep=3pt
\caption{The results of the new algorithm for Example \protect\ref
{exam6.2}}\label{table62}
\begin{tabular*}{\textwidth}{@{\extracolsep{\fill}}lccccc@{\hspace*{40pt}}c@{}}
\hline
& \multicolumn{1}{c}{\textbf{Part 1}} & \multicolumn{3}{c}{\textbf
{Part 2}}& & \\[-4pt]
& \multicolumn{1}{c}{\hrulefill} & \multicolumn{3}{c}{\hrulefill}&
& \\
\multicolumn{1}{@{}l}{$\bolds{j}$} & \multicolumn{1}{c}{$\bolds{\|\eps
_{j,1}\|^2}$} &\multicolumn{1}{c}{$\bolds{\|\eps_{j,2}\|^2}$} &
\multicolumn{1}{c}{$\bolds{T(\xi_j)}$}
& \multicolumn{1}{c}{$\bolds{\ds\frac{T (\xi_j)}{\|\eps_{j,2} \|^2}}$}
& \textbf{Support} & \textbf{Reference set}\\
\hline
$0$ & 1.25301\phantom{0} & & & & & $\cS\ot\{1,2,4,6,8,10\}$\\
$1$ & 0.040044 & & & & & \\
$2$ & 0.012839 & 0.008738 & 0.005404 & 0.6184 & $\{0.7,4.5,10\}$
& $\cS\ot\cS\cup\{0.7,4.5\}$ \\
$3$ & 0.006957 & 0.006827 & 0.006757 & 0.9897 & $\{0.5,3.6,10\}$
& $\cS\ot\cS\cup\{0.5,3.6\}$\\
$4$ & 0.006805 & 0.006797 & 0.006786 & 0.9996
& $\{0.5,3.4,10\}$ & $\cS\ot\cS\cup\{3.4\}$ \\
$5$ & 0.006793 & 0.006789 & 0.006786 & 0.9999 & & \\
$6$ & 0.006788 & 0.006787 & 0.006786 & 0.9999 & & \\
\hline
\end{tabular*}
\end{table}

%

\begin{exam}
\label{exam6.2}
In order to demonstrate that the algorithm can be used when dealing with
nonlinear regression models,\vadjust{\goodbreak} we consider two rival models
$ \eta_1(x,\theta) = \frac{\theta_{11}x}{x+\theta_{12}}$,
$\eta_2(x,\theta) = \theta_{21} ( 1 - e^{-\theta_{22}x})$,
where\vspace*{1pt} $\rho_{(1)} =(2.0, 1.0)$ and $\rho_{(2)}=(2.5, 0.5)$.
The weights in the criterion \eqref{det2} are $p_{1,2} = p_{2,1} = 1/2$.
The corresponding results are depicted in Table~\ref{table62} and the Newton method is started with
$\theta_{(1,2)} = (1,1)$, $\theta_{(2,1)} = (2,0.5)$
and $\cS=\{1,2,4,6,8,10\}$.
The degree of approximation is close to the optimum already
after 6 iteration steps, and the guaranteed efficiency is $99.9\%$.
The resulting design has masses $0.311 $, $ 0.415$ and $ 0.274$
at the points $ 0.5$, $3.4$ and $10.0$, respectively,
while the parameters of the solution of the approximation problem on
the right-hand side of \eqref{det5} are given (subject to rounding)
by the parameters
$\bar\theta_{(1,2)} = (3.008, 1.809)$,
and $\bar\theta_{(2,1)} = (1.721, 0.865)$.

\begin{figure}[b]

\includegraphics{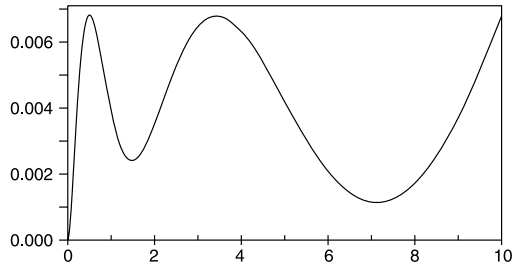}

\caption{The function $\psi$ in the equivalence theorem
(Theorem~\protect\ref{eqivthm})
for Example \protect\ref{exam6.2}.}\label{eqplot}
\end{figure}

The determination of the parameter $\theta^*$ that minimizes $T(\xi)$
as defined in \eqref{theta-xi} is done by Newton's method.
It yields the best $\theta$ in a neighborhood of the
computed solution.
Therefore, we have also performed an extensive global search for
the minimum and found a minimum that equals the result
of Newton's method up to rounding errors.
Now, the plot of the corresponding function $\psi$ in the
equivalence theorem (Theorem~\ref{eqivthm}) is shown
in Figure~\ref{eqplot}. We see that the design
is in fact $T_p$-optimal discriminating.
Note that the support of the resulting design
consists of $3$ points in accordance with Conjecture~\ref{con45}.
%
%
The corresponding results for the algorithm proposed by \citet
{atkfed1975a} are displayed in the
right part of Table~\ref{table61alt}. The algorithm needs $16$
iterations in order to find a design
with masses $ 0.32 $, $0.03$, $ 0.21 $, $0.12 $, $0.06$, $0.27$
at the (unnecessarily large set of) points
$ 0.5$, $3.0 $, $3.3$, $3.4 $, $ 3.8$, $ 10 $.
Here the lower bound of the efficiency is only $95.5\%$
if we take the best information from the previous steps.
{The new algorithm is obviously much faster.\vspace*{-2pt}}
\end{exam}

%
\begin{exam}
\label{exam73}
We consider $T_p$-optimal discriminating designs
for the four competing dose--response models
listed in Table~\ref{tabex1} in the \hyperref[sec1]{Introduction}
and the design space $\cX=[0,500]$.
Here, $d=6$ comparisons and $n=15$ parameters
are involved. Moreover, the model $\eta_4$ is nonlinear.
We use the weights $p_{i,j}=1/6$ if $i > j$ and $p_{i,j}=0$ otherwise
in the criterion \eqref{det2}.


The corresponding results are displayed in Table~\ref{table836}, which
shows that only $9$ iteration steps are required in order to obtain a
design with
at least $99.9 \%$ efficiency.
The resulting $T_p$-optimal discriminating design puts masses
$ 0.255$, $ 0.212$, $0.358$, $0.175$ at the points
$0 $, $78$, $245 $ and $500 $, respectively.

We finally note that we were not able to find a design with a
guaranteed efficiency of $80 \%$ using the algorithm proposed by \citet
{atkfed1975a}.\vspace*{-3pt}
%
\begin{table}
\def\arraystretch{0.9}
\tabcolsep=0pt
\caption{The results of the new algorithm for Example \protect\ref
{exam73}}\label{table836}
\begin{tabular*}{\textwidth}{@{\extracolsep{\fill}}lccccc@{\hspace*{40pt}}c@{}}
\hline
& \multicolumn{1}{c}{\textbf{Part 1}} & \multicolumn{3}{c}{\textbf
{Part 2}} & & \\[-4pt]
& \multicolumn{1}{c}{\hrulefill} & \multicolumn{3}{c}{\hrulefill} & &
\\
\multicolumn{1}{@{}l}{$\bolds{j}$} & \multicolumn{1}{c}{$\bolds{\|
\eps_{j,1}\|^2}$} &
\multicolumn{1}{c}{$\bolds{\|\eps_{j,2}\|^2}$} & \multicolumn
{1}{c}{$\bolds{T(\xi_j)}$}
& \multicolumn{1}{c}{$\bolds{\ds\frac{T (\xi_j)}{\|\eps_{j,2} \|^2}}$}
& \multicolumn{1}{c}{\textbf{Support}} & \multicolumn
{1}{c@{}}{\textbf{Reference set}}\\
\hline
$0$ & 16,661 & & & & &
$\cS\ot\{0,30,60,90,\ldots,450,500\}$\\

$1$ & 12,646 & & & & & \\
$2$ & \phantom{0.}9727 & 8923 & \phantom{0}275 & 0.0309 & $\{
0,50,290,450\}$
& $\cS\ot\cS\cup\{50,290\}$ \\
$3$ & \phantom{0.}8246 & 6901 & \phantom{0}764 & 0.1108 & $\{
0,60,290,450\}$ & \\
$4$ & \phantom{0.}5835 & 5081 & 2462 & 0.4846 & $\{0,70,260,500\}$
& $\cS\ot\cS\cup\{70,260\}$\\
$5$ & \phantom{0.}4543 & 4170 & 3016 & 0.7233 & $\{0,80,250,500\}$
& $\cS\ot\cS\cup\{80,250\}$\\
$6$ & \phantom{0.}4048 & 3619 & 3168 & 0.8754 & $\{0,80,240,500\}$
& $\cS\ot\{0,70,80,240,$\\
& & & & &
& $250,500\}$\\
$7$ & \phantom{0.}3446 & 3270 & 3194 & 0.9989 & & \\
$8$ & \phantom{0.}3201 & 3199 & 3195 & 0.9980 & $\{0,78,240,500\}$
& $\cS\ot\cS\cup\{78\}$ \\
$9$ & \phantom{0.}3197 & 3196 & 3195 & 0.9998 & & \\
\hline
\end{tabular*}\vspace*{-3pt}
\end{table}
\end{exam}

\section{Concluding remarks}

Our main theoretical result relates $T_p$-optimal discriminating
designs to
an approximation problem for vector-valued functions (Theorem~\ref{thm1}).
By duality theory we show that there exist
$T_p$-optimal designs with at most $n+1$ support points, where
$n$ is the number of parameters in the approximation problem (which
coincides with the total number of parameters of all regression
functions used in the comparisons).
These results are sufficient\vadjust{\goodbreak} if we are interested only one or two
comparisons among the rival models.
In this case the computations can be done
by an exchange-type algorithm
that was already proposed by \citet{atkfed1975a}.
This procedure is still the common tool
for dealing with design problems whenever $d=1$ or $d=2$.

The situation is different and the construction of
$T_p$-optimal discriminating designs becomes extremely difficult
and challenging if three or more comparisons are involved.
The number of support points can now be much smaller than $n+1$,
where $n$ is the total number of parameters of the models involved
in the $T_p$-optimality criterion.
Although a reduction of this number was already observed
in the case $d=2$, the amount of the reduction and its impact
become clear only when optimal discriminating design problems with $d
\ge3$
pairwise comparisons are studied.
For example, we have $n=15$ parameters in the dose-finding problems
listed in Table~\ref{tabex1}, but the support of the $T_p$-optimal discriminating design
consists of only 4 points.

Therefore, there are substantial differences between our new
algorithm and the generalization of the method by \citet{atkfed1975a}
beyond the case $d=1$.
Our algorithm is based on the related approximation problem
(Theorem~\ref{thm1}), and additionally we also add combinatorial aspects
[addition (iii) in Section~\ref{sec61}], which accelerate the speed of
convergence. Dual linear programs associated to small
subproblems determine the support of the resulting design
and prevent the algorithm from providing designs with too many support
points.
The masses are simultaneously computed by a stabilized version
of the equations in Corollary~\ref{weight},
while the commonly used algorithms
in each iteration step involve an update of the mass
at only one point and a renormalization.

Our numerical examples in Section~\ref{sec6} and in
the supplementary material [\citet{supp}] show that the new algorithm is able to
solve $T_p$-optimal discriminating design problems of higher
dimensions in situations
where the classical methods fail.

\section*{Acknowledgments}
We are very grateful to the referees and the Associate Editor for their
constructive comments on an earlier version of this manuscript.
In particular, one referee encouraged
us to include examples with a larger number of pairwise comparisons. By
these investigations we gained more insight in the optimization
problem, which led to a further improvement of the proposed algorithm.
We also want to thank Stefan Skowronek for providing
a code for the numerical calculations and Martina Stein, who typed
parts of this manuscript with considerable technical expertise.

\begin{supplement}[id=suppA]
\stitle{Optimal discriminating designs for several competing regression
models}
\slink[doi]{10.1214/13-AOS1103SUPP} 
\sdatatype{.pdf}
\sfilename{aos1103\_supp.pdf}
\sdescription{Technical details and more examples.\vadjust{\goodbreak}}
\end{supplement}

%
%

\printaddresses

\end{document}